\documentstyle[twoside,12pt]{article}

\newtheorem{thm}{\indent{\sc Theorem}}[section]
\newtheorem{defn}[thm]{\indent{\sc Definition}}

\newtheorem{prop}[thm]{\indent{\sc Proposition}} 
\newtheorem{cor}[thm]{\indent{\sc Corollary}}

\newcommand{\ssect}{\subsection}
\newcommand{\sssect}{\subsubsection}

\newcommand{\eqref}[1]{equation~(\ref{#1})}
\newcommand{\eref}[1]{\eqref{#1}}

\newcommand{~}{\nolinebreak[3] }

\newcommand{\seper}{\nopagebreak \begin{center}
	\underline{\hspace{2in}}
	\end{center}}

\newcommand{\bold}[1]{{\em (#1)}\index{#1}}
\newcommand{\bld}[1]{{\em (#1)}}

 \newcommand{\proof}[1]{\proo{#1}$\Box $}
\newcommand{\proo}[1]{{\em Proof: }#1}
\renewcommand{\Box}{ \hspace{1cm}{\rule{1.2ex}{2ex}}}

\newcommand{\RomCoeff}[2]{\romcoeff{#1}{#2}}
\newcommand{\rn}[1]{\left\lfloor #1 \right\rceil}

\newcommand{\romcoeff}[2]{\rn{{#1 \atop #2}}}

\newcommand{\eps}[1]{\rn{#1}_{\bullet}^{\epsilon }}
\newcommand{\ecoeff}[2]{\romcoeff{#1}{#2}_{\bullet}^{\epsilon }}

\newcommand{\ecoeffz}[2]{\romcoeff{#1}{#2}_{\bullet}^{0}}
\newcommand{\epso}[1]{\rn{#1}_{\bullet}^{1}}
\newcommand{\ecoeffo}[2]{\romcoeff{#1}{#2}_{\bullet}^{1}}

\newcommand{\fig}[2]{\begin{table}[htbp]
	\caption{#1} 
	\begin{center}\fbox{ \scriptsize \begin{tabular}#2 \end{tabular}}
	\end{center}
\end{table}}
\newcommand{\corn}{\multicolumn{1}{r|}{1}}
\newcommand{\ner}[1]{\cline{#1-#1}}
\newcommand{\figx}[2]{\begin{figure}[htbp] \caption{#1} 
	\begin{center} \scriptsize \mbox{}#2 \end{center} \end{figure}}

\begin{document}
\include{title}

\begin{abstract}
We pose the question of what is the best generalization of the factorial and
the binomial coefficient. We give several examples, derive their
combinatorial properties, and demonstrate their interrelationships.
\seper
{\bf G\'{e}n\'{e}ralisation des Coefficients du Bin\^ome}

On cherche ici \`{a} d\'{e}terminer est la meilleure
g\'{e}n\'{e}ralisation possible des factorielles et des coefficients
du bin\^ome. On s'interesse \`{a} plusieurs exemples, \`{a} leurs
propri\'{e}t\'{e}s combinatoires, et aux differentes relations qu'ils
mettent en jeu.
\end{abstract}

\begin{center}
{\it Dedicated to\\
David and Maureen}
\end{center}

\tableofcontents

\section{Introduction}

Despite being so fundamental to combinatorics, several authors have noticed
that ne is virtually unlimited in the choice of definition for the {\em
factorial}---at least as far as umbral calculus is concerned. Indeed, one is
presented with a bewildering number of alternatives each with its own notation.

We present a new definition of the factorial which generalizes the usual one, 
and study the binomial coefficients it induces. They are blessed with a variety
of combinatorial properties. However, what we are most interested is studying
the interrelationship between this factorial and other famous ones. 

 \ssect{The Roman Factorial} \index{Factorial,Roman}

We begin by presenting  a generalization of the
factorial $n!$ which makes sense for negative integral values of $n$
as well as nonnegative called the {\em  Roman factorial} $\rn{n}!$ after its
inventor Steve Roman. As usual for $n$ a nonnegative integer the factorial is
given by the product
$$ \rn{n}!=n!= 1\times 2 \times 3 \times \cdots \times n.$$
However, for $n$ a negative integer 
$$ \rn{n}!=  \frac{(-1)^{n+1}}{(-n-1)!} $$
\fig{Roman Factorials $\rn{n}!$}{{r|*{13}{c}}\index{Factorial,Roman}
\normalsize $n$ &\normalsize  $-6$ &\normalsize $-5$ &
\normalsize $-4$ &\normalsize $-3$ &\normalsize $-2$ & 
\normalsize $-1$& \normalsize 0 &\normalsize  1 &\normalsize  2 &
\normalsize  3 &\normalsize  4 &\normalsize  5 &\normalsize  6\\*
\hline  		
\normalsize ${\rn{n}}!$
& $-\frac{1}{120}$&$\frac{1}{24}$ &$- \frac{1}{6}$ & $\frac{1}{2}$
& $-1$ & 1 & 1 & 1 & 2 & 6 & 24 & 120 & 720}

{\begin{prop}[Knuth] \label{knuth}
For any integer $n$,
$$ \rn{n}!\rn{-n}! = (-1)^{n}|n|.\Box $$
\end{prop}}

More generally, for  every real number $a$, let
$$ \rn{a}!=\left\{\begin{array}{ll}
\Gamma (a+1)&\mbox{when  $a$ is not a negative integer, and}\\*[0.1in]
(-1)^{a-1}/(-a-1)!&\mbox{when $a$ is a negative integer}
\end{array} \right. $$
where $\Gamma (a)$ is the analytic Gamma function.

Thus, for all $a$
\begin{equation}\label{roma}
\rn{a}!/\rn{a-1}!=\rn{a}
\end{equation}
where {\em Roman } $a$ is defined to be
$$ \rn{a}=\left\{\begin{array}{ll}
a&\mbox{for $a\neq 0$}\\*[0.1in]
1&\mbox{for $a=0$.}
\end{array} \right. $$
Note that \eref{roma} and the condition $\rn{0}=1$ completely
characterizes the  Roman factorial of integers.

\ssect{The Roman Coefficients}
These extensions of the notion of factorial leads to a corresponding
generalization of the definition of Binomial Coefficients.

\begin{defn} \bold{Roman Coefficients} \label{roman}
For all real numbers $a$ and $b$, define
the {\em Roman coefficient} (read: ``Roman $a$ choose $b$'')
to be 
$$\RomCoeff{a}{b} = \frac{\rn{a}!}{\rn{b}! \rn{a-b}!}. $$
\end{defn}

\fig{Roman Coefficients, $\romcoeff{n}{k}$}{{r|rrrr|*{7}{r}}\label{RomanFig} 
\normalsize $n \backslash k$ &\normalsize $-4$ &\normalsize $-3$
&\normalsize  $-2$ &\normalsize $-1$ &\normalsize  0 &\normalsize  1
&\normalsize  2 &\normalsize  3 &\normalsize  4 &\normalsize  5
&\normalsize  6\\*
\hline 
\normalsize 6 & $-1/840$& 1/252 & $-1/56$ &1/7& 
1 & 6 & 15 & 20 & 15& 6&\corn\\* \ner{12}
\normalsize 5 &$-1/504$ &1/168& $-1/42$ &1/6 &
 1 & 5 & 10 & 10 & 5 &\corn &${1/6}$\\ \ner{11}
\normalsize 4 & $-1/280$ &1/105& $-1/30$ &1/5 & 1 & 4 & 6 & 4 & \corn &
${1/5}$& $-1/30$\\ \ner{10}
\normalsize 3& $-1/140$ &1/60& $-1/20$ &1/4& 1 & 3 & 3 & \corn & 
${1/4}$& $-1/20$ &1/60\\ \ner{9}
\normalsize 2 & $-1/80$ &1/30& $-1/12$ &1/3& 1 & 2 & \corn &
${1/3}$& $-1/12$ &1/30& $-1/60$\\ \ner{8}
\normalsize 1 & $-1/20$ &1/12& $-1/6$ &1/2&1 & \corn & ${1/2}$
& $-1/6$ &1/12& $-1/20$ &1/30\\ \ner{7}
\normalsize 0& $-1/4$ &1/3& $-1/2$ &1& \corn & ${1}$
 & $-1/2$ &1/3& $-1/4$ &1/5& $-1/6$\\ \hline 
\normalsize $-1$& $-1$ & 1 & $-1$& \corn & 1 & $-1$ & 1 & $-1$ & 1 & $-1$ &
1\\ \ner{5}
\normalsize $-2$ & 3 &$-2$ & \corn& ${-1}$
 & 1& $-2$ & 3 & $-4$ & 5 &$-6$&7 \\ \ner{4}
\normalsize $-3$& $-3$ & \corn & ${-1/2}$
 &$ -1/2$ & 1 & $-3$ & 6 & $-10$ & 15& $-21$ &28 \\ \ner{3}
\normalsize $-4$ & \corn & ${-1/3}$
&$-1/6$&$-1/3$&1 &$-4$ & 10 &$-20$ & 35& $-56$ &84 \\ \ner{2}
\normalsize $-5$ & ${-1/4}$ & 
$-1/12$ & $-1/12$ & $-1/4$ & 1 & $-5$ & 15 & $-35$ & 70&
$-126$ & 210}

When the two argument are both integers, the relationship between the Roman
coefficients and the binomial coefficients is given by the following:

{\begin{prop}\label{Rgns} \bld{The Six Regions}\index{Regions}
Let $n$ and $k$ be integers. Depending on what region
of the Cartesian plane the point $(n,k)$ is in, the
following formulas apply:

\begin{description}
\item[Region 1] If $n \geq k \geq 0$, then 
$\displaystyle \romcoeff{n}{k} = {n \choose
k}$.  
\fig{Region 1}{{r|*{8}{r}} 
\normalsize $n \backslash k$ &\normalsize  0 &\normalsize  1
&\normalsize  2 &\normalsize  3&\normalsize  4 &\normalsize  5
&\normalsize  6 &\normalsize  7\\ 
\hline 
\normalsize 7& 1 & 7 & 21 & 35 & 35 &21 & 7&1\\
\normalsize 6 & 1 & 6 & 15 & 20 & 15 & 6 & 1\\
\normalsize  5 & 1 & 5 & 10 & 10 & 5 & 1\\
\normalsize 4 & 1 & 4 & 6 & 4 & 1 \\
\normalsize 3& 1 & 3 & 3 & 1 \\
\normalsize 2 & 1& 2 & 1 \\
\normalsize  1 & 1 & 1 \\
\normalsize 0& 1}
\item[Region 2] If $k \geq 0 > n$, then 
$\displaystyle \romcoeff{n}{k} = 	(-1)^{k}
{{-n+k-1} \choose k}$.
\fig{Region 2}{{r|*{7}{r}}\label{Reg2}
\normalsize $n \backslash k$ &\normalsize  0 &\normalsize  1
&\normalsize  2 &\normalsize  3 &\normalsize  4 &\normalsize  5
&\normalsize  6\\ 
\hline 
\normalsize $-1$& 1 &$-1$ & 1 & $-1$ & 1 & $-1$ & 1\\
\normalsize $ -2$& 1 & $-2$ & 3 & $-4$ & 5 & $-6$&7 \\
\normalsize $-3$ & 1 & $-3$ & 6 & $-10$ & 15& $-21$ &28 \\
\normalsize $-4$ & 1 & $-4$ & 10& $-20$ & 35& $-56$ &84 \\
\normalsize $-5$ & 1 & $-5$ & 15 & $-35$& 70&$-126$ &210}
\item[Region 3] If $0 > n \geq k$, then 
$\displaystyle \romcoeff{n}{k} = 	(-1)^{n+k}
{{-k-1} \choose {n-k}}$.
\fig{Region 3}{{r|*{6}{r}}
\normalsize $n \backslash k$ &\normalsize $-6$&\normalsize
$-5$&\normalsize $-4$ &\normalsize $-3$ &\normalsize $-2$ &\normalsize $-1$\\ 
\hline 
\normalsize $-1$ &$-1$&1&$-1$ &1 &$-1$ & 1 \\
\normalsize $-2$&5& $-4$ & 3 & $-2$ & 1 \\
\normalsize $-3$ & $-10$ &6& $-3$ &1 \\
\normalsize $-4$&10& $-4$& 1\\
\normalsize $-5$&$-5$&1\\
\normalsize $-6$&1}
\item[Region 4] If $k > n \geq 0$, then 
$$\begin{array}{rcl}
\displaystyle \romcoeff{n}{k}&=& \displaystyle (-1)^{n+k} \frac{1}{n-k} {k
\choose n} ^{-1} 
= \displaystyle (-1)^{n+k+1}\frac{1}{n+1}{k\choose n+1 }^{-1}\\*
&=& \displaystyle (-1)^{n+k+1}\frac{1}{k}{k-1 \choose n }^{-1}\\*
&=& \displaystyle (-1)^{n+k+1}n\sum _{j\geq 0}S(j,n)/k^{j+1}
= \displaystyle (-1)^{n+k}\left[\Delta ^{n}\frac{1}{x-k}\right]_{x=0}
\end{array}$$
where the $S(j,n)$ are the Stirling numbers of the second kind, and $\Delta $
is the forward difference operator $\Delta p(x)=p(x+1)-p(x)$.
\fig{Region 4}{{r|*{7}{c}}
\normalsize $n \backslash k$ &\normalsize  1 &\normalsize  2
&\normalsize  3 &\normalsize 4 &\normalsize  5 &\normalsize  6
&\normalsize  7\\
\hline 
\normalsize 6 & &&&&&& 1/7\\
\normalsize 5 & &&&&&1/6 &$ -1/42$\\ 
\normalsize 4 & &&&&1/5&$-1/30$&1/105\\
\normalsize 3&&&&1/4&$-1/20$ &1/60& $-1/140$\\
\normalsize 2 &&&1/3 & $-1/12$ &1/30& $-1/60$ & 1/105\\ 
\normalsize 1 && 1/2 & $-1/6$ &1/12& $-1/20$ &1/30&$-1/42$\\ 	
\normalsize 0 &1 & $-1/2$ &1/3 & $-1/4$ &1/5& $-1/6$ &1/7}
\item[Region 5] If $n \geq 0 > k$, then 
$$\begin{array}{rcl}
\displaystyle \romcoeff{n}{k} & =&
\displaystyle (-1)^{k}\frac{1}{k} {n-k \choose n} ^{-1}
= \displaystyle (-1)^{k}\frac{1}{k-n}{n-k-1\choose n }^{-1}\\*
&=&\displaystyle (-1)^{k+1}\frac{1}{n+1}{n-k-1 \choose n+1}^{-1}\\
&=&\displaystyle \romcoeff{n}{n-k} = (-1)^{k}\left[\Delta
^{n}\frac{1}{x-n+k}\right]_{x=0}\\*
&=&\displaystyle  -B(k-n,-k)
\end{array}$$
where the pair $(n,n-k)$ lies in region 4 (defined above), and $B(n,k)$ is the
analytic Beta function.
\fig{Region 5}{{r|*{4}{l}}
\normalsize $n \backslash k$ &\normalsize  $-4$ &\normalsize$-3$
&\normalsize  $-2$ &\normalsize  $-1$ \\
\hline 
\normalsize 6&$-1/840$&1/252&$-1/56$&1/7\\
\normalsize 5&$-1/504$&1/168&$-1/42$&1/6\\
\normalsize 4&$-1/280$&1/105&$-1/30$&1/5\\
\normalsize 3&$-1/140$&1/60&$-1/20$&1/4\\
\normalsize 2&$-1/80$&1/30&$-1/12$&1/3\\
\normalsize 1&$-1/20$&1/12&$-1/6$&1/2\\
\normalsize 0&$-1/4$&1/3&$-1/2$&1}
\item[Region 6] Region 6: If $0>k>n$, then 
$$\begin{array}{rcl}
\displaystyle \romcoeff{n}{k} 
&=&\displaystyle \frac{1}{n-k} {-n-1\choose -k-1} ^{-1} 
=\displaystyle \frac{1}{k}{-n-1 \choose -k }^{-1}\\*
&=&\displaystyle \frac{1}{n+1}{-n-2 \choose -k-1 }^{-1}
=\displaystyle \left[\Delta ^{k-n-1}\frac{1}{x+n+1}\right]_{x=0}\\*
&=&\displaystyle (-1)^{k+1}\romcoeff{k-n-1}{-n-1} 
=\displaystyle (-1)^{k+1}\romcoeff{k-n-1}{k} 
\end{array}$$
where the pair $(k-n-1,-n-1)$ lies in region 4 (defined above), and the pair
$(k-n-1,k)$ lies in 
region 5 (defined above).$\Box $
\fig{Region 6}{{r|*{6}{l}}
\normalsize $n \backslash k$ &\normalsize  $-6$ &\normalsize $-5$
&\normalsize  $-4$ &\normalsize $-3$ &\normalsize  $-2$ &\normalsize  $-1$\\ 
\hline 
\normalsize $-2$&&&&&& $-1$ \\
\normalsize $-3$&&&&& $-1/2 $&$ -1/2$\\
\normalsize $-4$&&&&$-1/3$ & $-1/6$ &$ -1/3$\\
\normalsize$ -5$&&&$-1/4$ & $-1/12$& $ -1/12$ & $-1/4$\\
\normalsize $ -6$&&$-1/5$&$-1/20$&$-1/30$&$-1/20$&$-1/5$\\
\normalsize $-7$&$-1/6$&$-1/30$&$-1/60$&$-1/60$&$-1/30$&$-1/6$}
\end{description}
\end{prop}

{Note that in regions 1, 2, and 3, the Roman coefficients
equal binomial
coefficients up to a permutation and a
 change of sign. 
In regions 4, 5, and 6, the Roman coefficients are expressed simply in terms of
the {\it reciprocals} of the binomial coefficients. Furthermore,
regions 4, 5, and 6 are identical up to permutation and change of
sign. Thus, all of the Roman coefficients are related in a
simple way to those in the first quadrant (regions 1 and 4).
In particular, the Roman coefficients always equal
 integers or the reciprocals of integers.}

\ssect{Properties of Roman Coefficients}\label{props}
\index{Roman Coefficients}

Several binomial coefficient identities extend to Roman
coefficients.  

\begin{prop} \bold{Complementation Rule} 
For all  real numbers $a$ and $b$, 
$$ \RomCoeff{a}{b} = \RomCoeff{a}{a-b}.\Box $$ 
\end{prop} 

\begin{prop} \bold{Iterative Rule}\label{iterated}  
For all real numbers $a$, $b$, and $c$; 
$$ \RomCoeff{a}{b} \RomCoeff{b}{c}
= \RomCoeff{a}{c} \RomCoeff{a-c}{b-c}.\Box $$ 
\end{prop}

\begin{prop}\bold{Pascal's Recursion}\label{Pascal} 
If $a$ and $b$ are distinct and nonzero real numbers, then we have
\[ \RomCoeff{a}{k} =
\RomCoeff{a-1}{k} + \RomCoeff{a-1}{k-1}.\] 
\end{prop}

\proo{Since under these
conditions $\rn{a}=a$, $\rn{b}=b$, 
and $\rn{a-b}=a-b$,
\begin{eqnarray*}
\RomCoeff{a-1}{b}+\RomCoeff{a-1}{b-1} &=&
\frac{\rn{a-1}!}{\rn{a-b-1}!\rn{b}!} + \frac{\rn{a-1}!}{\rn{a-b}!\rn{b-1}!}\\*
&=& \rn{a-b}\left(\frac{\rn{a-1}!}{\rn{a-b}!\rn{b}!}\right)
+\rn{b}\left(\frac{\rn{a-1}!}{\rn{a-b}!\rn{b}!}\right) \\
&=& \rn{a}\left(\frac{\rn{a-1}!}{\rn{a-b}!\rn{b}!}\right)\\
&=&\frac{\rn{a}!}{\rn{a-b}!\rn{b}!}\\*
&=& \RomCoeff{a}{b}.\Box 
\end{eqnarray*}}

{\begin{cor} 
If  $r$ is a nonnegative integer  and the pairs of integers 
$(n,k)$, $(n+r,k)$, $(n,k+1)$, and $(n+r+1,k+1)$ all lie in the 
same region (as defined in Theorem~\ref{Rgns}), we have
$$ \sum _{m=n}^{n+r}\romcoeff{m}{k}=\romcoeff{n+r+1}{k+1}
-\romcoeff{n}{k+1}. $$
\end{cor}}

\proof{Induction on $r$.}

Contrast this corollary with this classical result involving binomial
coefficients in which for $n\geq k\geq 0$,
$$\sum _{m=k}^{n}{m \choose k }={n+1 \choose k+1 }.  $$

Analogous results hold more generally for real numbers.

If we adopt Iverson's\index{Iverson} notation for the moment writing
logical expressions in parenthesis to mean 1 if true and 0
if false, then in the discrete case we have the following proposition:

{\begin{prop}\label{rotref}\index{Knuth{,} Donald,
Rotation/Reflection Law}
\bld{Knuth's Rotation/Reflection Law} For any integers $n$
and $k$,
$$ (-1)^{k+(k>0)}\romcoeff{-n}{k-1} =
(-1)^{n+(n>0)}\romcoeff{-k}{n-1}.  $$
\end{prop}}

\proo{By Proposition~\ref{knuth}, we have 
$$ \romcoeff{n}{k}=
(-1)^{n+k+(n<0)+(k<0)}\romcoeff{-k-1}{-n-1}.\Box  $$}

{\begin{prop}\bold{Roman's Identity}
For all integers $n$ and $k$,
$$ \romcoeff{n}{k}\romcoeff{k}{n} = \frac{(-1)^{n+k}}{|n-k|}. $$
\end{prop}}

\proof{Proposition~\ref{knuth}.}

\ssect{Generalizations of the Roman Coefficients}
\index{Roman Coefficients}
\label{Gen}

The Roman coefficients defined earlier were very useful. However, there are
several other generalizations of binomial 
coefficients. For example, recall the classical definition of extended
binomial coefficients.

\begin{defn}\bld{Classical Extended Binomial
Coefficient}\index{Binomial Coefficient}
 Given a field element $x\in K$ in a field $K$ of
characteristic zero, and a nonnegative integer $k$, 
define the binomial coefficient ``x choose k'' to be: \[ {x
\choose k} = (x)_k/k!.\] 
where $(x)_{k}$ denotes the lower factorial of $x$ of degree $k$
$$(x)_k = \left\{ \begin{array}{rcll}
\displaystyle \prod_{i=0}^{k-1}(x-i) & = & x(x-1)\cdots
(x-k+1) & \mbox{for $k \geq 0$, and} \\
\\
\displaystyle \prod _{i=k}^{-1}(x-i)^{-1} & = &
{1}/{(x+1)(x+2)\cdots (x-k)}
& \mbox{for $k < 0$.} \end{array} \right.$$ 
\end{defn}

What is the relationship between the Roman coefficients and the other
generalizations of binomial coefficients?
To fully answer this
question, we must generalize our notion of harmonic factorial.

\sssect{Knuth Coefficients}\label{knsec}

Adopt the following convention independently discovered by
Donald Knuth.\index{Knuth{,} Donald} 

\begin{defn}\bld{Knuth
Factorial}\label{epsilon}\index{Factorial,Knuth}\index{Knuth{,} Donald,
Factorial} 
Define $\eps{a}$ for $a$ a real number to be the most
significant term of $\Gamma (a+1+\epsilon )$ where $\epsilon $
is an infinitesimal in from the field of surreal
numbers (a non-Euclidean
field which contains the real numbers). 
\end{defn}

 Thus, for $a$ real,
\begin{equation}\label{KF}
\eps{a}=\left\{ \begin{array}{ll}
\Gamma (a+1)&\mbox{when $a$ is not a negative integer, and}\\[0.1in]
(-1)^{a-1}\omega /(-a-1)!&\mbox{when $a$ is a negative integer.}
\end{array} \right.
\end{equation}
where $\omega =1/\epsilon $. This choice of factorial whoold have led to
``tags'' of $\epsilon $ or $\omega $ in appropriate places
in results of this paper.

For instance, again for $a$ real, 
$$ \rn{a}^{\epsilon }=\left\{\begin{array}{ll}
a&\mbox{if $a\neq 0$, and }\\[0.1in]
\epsilon &\mbox{if $a=0$.}
\end{array} \right. $$
This is perhaps more natural since then $\rn{a}^{\epsilon}$ only
differs from $a$ by at most an infinitesimal.

If we adopt \eref{KF} as our definition where
$\epsilon $ can be any arbitrary constant, then
the Roman factorial  can be seen as a 
special case of the Knuth factorial where $\epsilon =1$.
That is, $\rn{a}!=\epso{a}$. Thus, the
motivation for our notation.

Let us proceed to generalize the Roman coefficients.

\begin{defn}\bld{Knuth Coefficient}\index{Knuth{,} Donald, Coefficient}
For all $a$ and $b$, define the {\em Knuth coefficient}
$\ecoeff{a}{b}$ by the fraction
$$\ecoeff{a}{b}=\frac{\eps{a}}{\eps{b}\eps{a-b}}.$$
\end{defn}

Clearly, $\ecoeffo{a}{b}=\romcoeff{a}{b}$. 

Let us
calculate $\ecoeff{n}{k}$ for each of the six regions mentioned in
Theorem~\ref{Rgns}.

\begin{prop}\bld{The Six Regions}\index{Regions}
Let $\epsilon $ be a nonzero complex number or surreal
number, and and $n,k$ be integers. Depending
on what region of the Cartesian plane the pair $(n,k)$ is in, the following
formulas apply:
\begin{description}
\item[Region 1]  If $n\geq k\geq 0$, then $$\displaystyle \ecoeff{n}{k}
=\displaystyle {n \choose k}=\displaystyle \romcoeff{n}{k}.$$
\item[Region 2]  If $k\geq 0>n$, then
$$ \ecoeff{n}{k} =(-1)^{k}{-n+k-1\choose k } = \romcoeff{n}{k}.$$
\item[Region 3]  If $0>n\geq k$, then 
$$ \displaystyle \ecoeff{n}{k} =  \displaystyle (-1)^{n+k}{-k-1\choose n-k } =
\displaystyle \romcoeff{n}{k}.$$
\item[Region 4] If $\displaystyle k>n\geq 0$, then
$$\displaystyle \ecoeff{n}{k} =\displaystyle \frac{(-1)^{n+k}\epsilon
}{(n-k){k\choose n }} =\displaystyle \epsilon \romcoeff{n}{k}.$$
\item[Region 5]  If $n\geq 0>k$, then 
$$ \displaystyle \ecoeff{n}{k}
=\displaystyle \frac{(-1)^{n+k}\epsilon }{k{n-k\choose n }}
=\displaystyle \epsilon \romcoeff{n}{k}. $$ 
\item[Region 6]  If $0>k>n$, then 
$$\displaystyle \ecoeff{n}{k}  = \displaystyle \frac{\epsilon}{k{-n-1\choose -k
}} =\displaystyle \epsilon \romcoeff{n}{k}.\Box$$
\end{description}
\end{prop}

\sssect{Gamma-Coefficients}
A limiting case of the Knuth coefficient is of special interest.
\begin{defn}\bold{Gamma-Coefficient}
Let $n$ and $k$ be arbitrary integers.
Define the Gamma-Coefficient 
$$\ecoeffz{n}{k}= \lim_{\epsilon \rightarrow 0}\ecoeff{n}{k} = \lim_{\epsilon
\rightarrow 0} \frac{\Gamma (n+1+\epsilon )}{\Gamma (k+1+\epsilon )\Gamma
(n-k+1+\epsilon )}. $$ 
\end{defn}

Note however that $\ecoeff{-1}{1/2}$ diverges as $\epsilon $
tends to zero, so it is impossible to define a
Gamma-Coefficient $\ecoeffz{a}{b}$ for $a$ and $b$ real.

\fig{Gamma-Coefficient $\ecoeffz{n}{k}$}{{r|rrrr|*{7}{r}}\label{GammaFig}
\normalsize $n \backslash k$ &\normalsize  -4 &\normalsize  -3
&\normalsize  -2 &\normalsize  -1 &\normalsize  0 &\normalsize  1
&\normalsize  2 &\normalsize  3 &\normalsize  4 &\normalsize
\normalsize  5 &\normalsize  6\\
\hline  
\normalsize 6 & 0&0&0&0 & 1 & 6 & 15 & 20 & 15 & 6 & \corn\\ \ner{12}
\normalsize 5 & 0&0&0&0& 1 & 5 & 10 & 10 & 5 & \corn &0\\ \ner{11}
\normalsize 4 & 0&0&0&0& 1 & 4 & 6 & 4 & \corn &0&0\\ \ner{10}
\normalsize 3 & 0&0&0&0& 1 & 3 & 3 & \corn &0&0&0\\ \ner{9}
\normalsize 2 & 0&0&0&0& 1 & 2 & \corn &0&0&0&0\\ \ner{8}
\normalsize 1 & 0&0&0&0& 1 & \corn & 0&0&0&0&0\\ \ner{7}
\normalsize 0 & 0&0&0&0& \corn & 0&0&0&0&0&0\\ \hline 
\normalsize $-1$& $-1$ & 1 & $-1$ & \corn & 1 & $-1$ & 1 & $-1$ & 1 & $-1$ &
1\\ \ner{5} 
\normalsize $-2$& 3 & $-2$ & \corn & 0 & 1 & $-2$ & 3 & $-4$ & 5 &$-6$&7 \\
\ner{4} 
\normalsize $-3$& $-3$ & \corn & 0&0& 1 & $-3$ & 6 & $-10$ & 15&$-21$&28 \\ \ner{3}
\normalsize $-4$& \corn & 0&0&0& 1 &$-4$ & 10 & $-20$ & 35&$-56$&84 \\ \ner{2}
\normalsize $-5$& 0&0&0&0& 1 & $-5$ & 15 & $-35$ & 70&$-126$&210} 

In regions 1, 2, and 3, the Gamma-Coefficients are equal to the Roman
coefficients. In regions 4, 5, and 6, the Gamma-Coefficients
are identically zero whereas the Roman coefficients are never zero.
Nevertheless, one should note that even when the Classical
binomial coefficient and the Roman coefficient differ, the difference
is at most one. 

Also, notice that the Gamma-Coefficients are always integers. In
particular, for $k\geq 0$ ({\it i.e.:} regions 1, 2, and 3), the
Gamma-Coefficients agree with the classical extended binomial coefficients.

The identities mentioned in \S\ref{props} generalize to
Gamma-Coefficients.  However, we defer any discussion of the combinatorial
significance to \cite{hybrid}.

\begin{prop}\bold{Complementation Rule}
For all real numbers $a$, $b$, and $\epsilon$,
$\ecoeff{a}{b}=\ecoeff{a}{a-b}$. In particular, for all integers $n$ and $k$,
$\ecoeffz{n}{k}=\ecoeffz{n}{n-k}.\Box $
\end{prop}

\begin{prop}\bold{Iterative Rule}\label{gammaprod}
For all real numbers $a,$ $b,$ $c,$ and $\epsilon $,
$$\ecoeff{a}{b}\ecoeff{b}{c}=\ecoeff{a}{c}\ecoeff{a-c}{b-c}.$$
In particular, for all integers $m$, $n$, and $k$,
$$\ecoeffz{m}{n} \ecoeff{n}{k}= \ecoeffz{m}{k}
\ecoeffz{m-k}{n-k}.\Box $$
\end{prop}
 
\begin{prop}\bold{Pascal's Recursion}
\begin{enumerate}
\item  Let $a$ and $b$ be distinct nonzero real numbers, and let
$\epsilon $ be a nonzero complex number.
Then
$$ \ecoeff{a}{b}=\ecoeff{a-1}{b}+\ecoeff{a-1}{b-1}.$$
\item For all $n$ and $k$,
$$ \ecoeffz{n}{k}=\ecoeffz{n-1}{k}+\ecoeffz{n-1}{k-1}$$
unless  $n=k=0.\Box$ 
\end{enumerate}
\end{prop}

Nevertheless, $\ecoeffz{0}{0}=1$ whereas
$\ecoeffz{-1}{-1}+\ecoeffz{-1}{0}=1+1=2$.

\sssect{Other Factorials}\index{Factorials,Other}

Actually as noted by Ueno \cite{Ueno} and Roman \cite{RomanB,Rom2}, any choice
of $\rn{a}!$ could be used for computations involving an umbral calculus. The
only restrictions are that 
$\rn{0}!$ must equal one, and for the so called continuous iterated logarithmic
algebra of \cite{ch4},  the function $a\mapsto \rn{a}!$ must be
continuous. 

For example, if we chose $\rn{a}!=1$ as in \cite{DiB}, then we have the theory
of convolution sequences. 

Whereas, if for $n$ an integer, we set as in \cite{More}\label{qref}
$$ \rn{\rn{n}}= \frac{q^{\rn{n}}-1}{q-1}, $$
then we achieve a $q$-analog of the ``$\rn{n}$-Logarithmic
theory.''

\sssect{Multinomial Coefficients}

Recall the usual definition of a multinomial coefficient.

\begin{defn}\bold{Classical Multinomial Coefficient}\label{multi}
Let  $n$ be a nonnegative integer, and let $\beta $ be a
vector with finite support of nonnegative integers. Then
define the multinomial coefficient 
$n$ choose $\beta $ to be 
$$ {n\choose \beta }=\left\{\begin{array}{ll}
\displaystyle n!\left(\prod _{k} \beta _{k}! \right)^{-1}&
\mbox{if }\displaystyle |\beta |=n \mbox{, and}\\[0.1in]
0&\mbox{otherwise.}
\end{array}\right. $$
\end{defn}

Note that ${n \choose \beta }$ is the number of ordered
partitions of type $\beta $ of a given $n$-set.

By analogy, for all reals $a$, and all real vectors $\beta $ with
finite support, define the {\em multinomial Roman
coefficient} $a$ choose $\beta $ to be
$$ \RomCoeff{a}{\beta } = 
\left\{\begin{array}{ll}
\displaystyle \rn{a}! \left(\prod_{k} \rn{\beta _{k}}!\right)^{-1}&
\mbox{if }\displaystyle |\beta |=a\mbox{, and}\\[0.1in]
0&\mbox{otherwise,}
\end{array}\right.$$
Define the multinomial Knuth coefficients\index{Knuth{,} Donald, Coefficient}
and Gamma-Coefficients\index{Gamma-Coefficient} similarly.
The multinomial Gamma-coefficients are well defined since they would
only diverge if some denominator had an excess of factors
of $\epsilon$.
However, that could only happen if $n<0$ and $k_{i}\geq 0$ for all
$i$, but in that case, $n\neq \sum _{i=1}^{j}k_{i}$, so the multinomial
$\epsilon $-coefficient, $\ecoeff{n}{(k_{i})_{i=1}^{j}}$ is zero by
definition. Contradiction! Thus, the Gamma-coefficients are well defined.

In terms of multinomial coefficients, Proposition~\ref{iterated} becomes 
$$ \romcoeff{n}{k}\romcoeff{k}{r}=\romcoeff{n}{n-k,k-r,r}
=\romcoeff{n}{r}\romcoeff{n-r}{k-r},$$ 
and Proposition~\ref{gammaprod}
becomes 
$$ \ecoeffz{n}{k}\ecoeffz{k}{r}=\ecoeffz{n}{n-k,k-r,r}=\ecoeffz{n}{r}
\ecoeffz{n-r}{k-r}, $$ 
More generally, we have the
following theorem.
\begin{prop}\bold{Iterative Rule}
Let $(k_{i})_{i=1}^{j}$ be a finite sequence
of integers with sum $n$. Then\begin{eqnarray*}
\romcoeff{n}{(k_{i})_{i=1}^{j} }&=&\prod _{m=2}^{j} \romcoeff{\sum
_{i=1}^{m}k_{i}}{k_{m}},\\*
\ecoeff{n}{(k_{i})_{i=1}^{j} }&=&\prod _{m=2}^{j} \ecoeff{\sum
_{i=1}^{m}k_{i}}{k_{m}}\mbox{, and}\\*
\ecoeffz{n}{(k_{i})_{i=1}^{j} }&=&\prod _{m=2}^{j} \ecoeffz{\sum
_{i=1}^{m}k_{i}}{k_{m}}.\Box
\end{eqnarray*} 
\end{prop}

As opposed to ordinary Roman coefficients, these multinomial Roman
coefficients are not always integers or reciprocals of integers---even when all
of the arguments are integers. For example, 
$\romcoeff{3}{2,2,-1} = \frac{3}{2}$.

However,  the multinomial Gamma coefficients are 
always integers, for if $\ecoeffz{n}{(k_{i})_{i=1}^{j}}$ is nonzero, then
we are in one of the following two cases. Either $n\geq 0$, and
$k_{i}\geq 0$ for all $i$, or $n<0$ and there is a unique $i$ such
that $k_{i}<0$. In the first case, these are ordinary multinomial
coefficients. It  suffices to consider the other case. Thus, $n<0$.
Without loss of generality, let $k_{1}<0$. Now, 
$$\ecoeffz{n}{(k_{i})_{i=1}^{j}}=
(-1)^{n+k_{1}}\ecoeffz{-k_{1}-1}{-n-1,k_{2},\cdots,k_{j}}
= (-1)^{n+k_{1}}{-k_{1}-1 \choose -n-1,k_{2},\cdots,k_{j}}$$
where $-k_{1}-1,-n-1,k_{2},k_{3},\cdots,k_{j-1},k_{j}\geq 0$.
Hence, all the nonzero multinomial Gamma coefficients are (up to sign)
ordinary multinomial coefficients, and thus integers.

\ssect{Resistance of the $n$-cube}
\index{Cube}\index{Electronics}\index{Resistance}\index{Combinatorial
Interpretations} 
Via the Gamma-coefficients and the theory of sets with a negative number of
elements \cite{hybrid}, we have a simple combinatorial interpretations for the
Roman coefficients in regions 1, 2, and 3. However, what is the significance of
the Roman coefficients in regions 4, 5, and 6? In these regions, the Roman
coefficients are the reciprocals of integers, so they do not enumerate any set.
However, the following application illustrates their  combinatorial
significance.  

\begin{prop}\label{cube}
Consider an $n$-cube
\figx{The $n$-cube}{\vspace{2in}}
in which each edge is
represented by a wire of resistance 1$\Omega $ (one Ohm). The
resistance between two opposing vertices of the cube is 
$$R_{n}=2^{-n}\sum _{i=1}^{n}i^{-1}2^{i}\Omega.$$
\end{prop}
\fig{Resistance of the $n$-cube}{{r|*{8}{c}}
\normalsize $n$ & \normalsize 0 &\normalsize  1 &\normalsize  2 &
\normalsize  3 &\normalsize  4 &\normalsize  5 &\normalsize  6 & \normalsize
7\\  
\hline \normalsize
${R_{n}}$&0&1&1&${5}/{6}$ & ${2}/{3} $ &
${8}/{15}$ &${13}/{30}$ & $151/340$}

\proof{The cube is isomorphic to the Hasse diagram of the boolean lattice of
subsets of  $\{1,2,\ldots ,n \}$. Without loss of
generality, the two opposing vertices are $\emptyset$, and
$\{1,2,\ldots ,n \}$. To compute the resistance, connect these two vertices
to a 1V battery. The resulting current (in Amperes) is equal to the
resistance (in Ohms).

By symmetry, each vertex on level $i$ of the lattice has the same
potential. Hence, we can consider each level as a single node without
effecting the resistance. Any two adjacent levels $i$ and $i+1$ are
connected by $(n-i){n\choose i }$ edges. Thus, the resistance between
levels $i$ and $i+1$ is $\frac{1}{n-i}{n\choose i }^{-1}\Omega $, or in
the notation of Roman coefficients, the resistance between levels $i$
and $i+1 $ is $(-1)^{n+i}\romcoeff{i}{n}\Omega
=-\romcoeff{-n-1}{1-i}\Omega $. 

The total resistance $R_{n}$ is the sum of the resistances between the
adjacent levels, 
$$ R_{n}= -\sum _{i=-n}^{-1}\romcoeff{-n-1}{i}\Omega .$$ 
By Theorem~\ref{Pascal},
$$ 2R_{n}=R_{n-1}+\frac{2}{n}\Omega . $$ 
We conclude by  induction noting
that $R_{0}=0$ and $R_{1}=1\Omega $.}

Note that as $n$ tends towards infinity, $R_{n}$ tends towards zero as $2/n$.

\end{document}